\documentclass[11pt, reqno]{amsart}

\usepackage{amsmath}
\usepackage{amscd}
\usepackage{mathrsfs}
\usepackage{amsfonts}
\usepackage{amssymb}
\usepackage{enumerate}
\usepackage{setspace}

\usepackage{array}

\usepackage{color}

\usepackage[hyperindex, colorlinks=true]{hyperref}

\usepackage[margin=1.2in]{geometry}

\newcommand{\eqn}{\begin{eqnarray}}
\newcommand{\een}{\end{eqnarray}}

\newcommand{\pat}{\partial_t}

\newtheorem{theorem}{Theorem}[section]
\newtheorem{corollary}[theorem]{Corollary}
\newtheorem{lemma}{Lemma}[section]

\theoremstyle{definition}
\newtheorem{definition}{Definition}[section]

\newtheorem{remark}{Remark}

\numberwithin{equation}{section}

\begin{document}

\title{Global existence for some transport equations with nonlocal velocity}

\author{Hantaek Bae}
\address{Department of Mathematics, University of California, Davis, USA}
\email{hantaek@math.ucdavis.edu}

\author{Rafael Granero-Belinch\'{o}n}
\address{Department of Mathematics, University of California, Davis, USA}
\email{rgranero@math.ucdavis.edu}

\dedicatory{This paper is dedicated to Professor Eitan Tadmor on the occasion of his 60th birthday.}

\date{\today}

\begin{abstract}
In this paper, we study transport equations with nonlocal velocity fields with rough initial data.  We address the global existence of weak solutions of an one dimensional model of the surface quasi-geostrophic equation and the incompressible porous media equation, and one dimensional and $n$ dimensional models of the dissipative quasi-geostrophic equations and the dissipative incompressible porous media equation.  
\end{abstract}

\maketitle

\section{Introduction} \label{sec:1}
In this paper, we study several active scalar equations with nonlocal velocity fields. Here, the non-locality means that the velocity field is defined through a singular integral operator that is represented in terms of  a Fourier multiplier. For example, in the two dimensional Euler equation in vorticity form  \cite{Bertozzi}, the velocity is recovered from the vorticity $\omega$ through
\eqn \label{2DEuler}
u=\nabla^\perp(-\Delta)^{-1}\omega \quad \text{or equivalently} \quad  \widehat{u}(\xi)=\frac{i\xi^{\perp}}{|\xi|^{2}}\widehat{\omega}(\xi).
\een
Other nonlocal and quadratically nonlinear equations appear in many applications. Prototypical examples are the surface quasi-geostrophic equation, the incompressible porous medium equation, Stokes equations, magnetogeostrophic equation and their variants.  We briefly introduce the equations below.

\subsubsection*{The surface quasi-geostrophic equations} 
The surface quasi-geostrophic equation describes the dynamics of the mixture of cold and hot air and the fronts between them in 2 dimensions \cite{Constantin2, Pedlosky}. The equation is of the form
\begin{equation}\label{SQG}
\theta_{t}+u\cdot\nabla\theta=0, \quad u=\left(-\mathcal{R}_{2}\theta, \mathcal{R}_{1}\theta\right),
\end{equation} 
where the scalar function $\theta$ is the potential temperature and $\mathcal{R}_{j}$ is the Riesz transform
\[
\mathcal{R}_{j}f(x)=\frac{1}{2\pi} \text{p.v.} \int_{\mathbb{R}^{2}} \frac{(x_{j}-y_{j})f(y)}{|x-y|^{3}}dy, \quad j=1,2. 
\]
Similar model equations of \eqref{SQG} with different types of nonlocal velocities are proposed and analyzed in \cite{Balodis},  \cite{Caffarelli} and \cite{Chae2singularSQG}, respectively (see also \cite{Berselli, CaffarelliVasseur, Chaeeqantonio, KiselevSQG, Li, Omar, Omar2}): 
\begin{equation*}
 \begin{split}
 & \theta_{t}+u\cdot\nabla\theta=0, \quad u=\mathcal{R}\theta,\\
 & \theta_{t} +\nabla \cdot \left(\theta \mathcal{R}\theta \right)=0,\\
 & \theta_{t}+u\cdot\nabla\theta=0, \quad u=\nabla^\perp\Lambda^{\beta-2}\theta,\quad 1<\beta\leq2.
 \end{split}
\end{equation*}
We finally introduce the two dimensional Euler-$\alpha$ model in vorticity form: 
\eqn \label{Euleralpha}
\theta_{t}+u\cdot\nabla\theta=0, \quad  u=\nabla^\perp \Lambda^{-2+\alpha}\theta, \quad \alpha \in[0,1]
\een 
which interpolates between \eqref{2DEuler} ($\alpha=0$) and \eqref{SQG} ($\alpha=1$).

\subsubsection*{The incompressible porous medium equation}
This equation takes the form \cite{castro2009incompressible}
\eqn \label{IPM}
\theta_{t}+u\cdot\nabla\theta=0, \quad u=\mathcal{R}^{\perp}\mathcal{R}_1\theta,
\een 
where $\theta$ is now the density of the incompressible fluid moving through a homogeneous porous medium. The following version also has been studied with $\beta>0$ \cite{Friedlander}:
\eqn \label{IPMsingular}
\theta_{t}+u\cdot\nabla\theta=0, \quad u=\Lambda^\beta\mathcal{R}^{\perp}\mathcal{R}_1\theta.
\een 
Another important equation related to flow in porous media is the Stokes equation: 
\eqn \label{Stokes}
\theta_{t}+u\cdot\nabla\theta=0, \quad -\Delta u=-\nabla p-(0,\theta)^t,\quad \nabla\cdot u=0,
\een
where the velocity field is given by 
\begin{equation*} 
 \begin{split}
 \text{2D}: \ & u=(-\Delta)^{-1}\mathcal{R}^{\perp}\mathcal{R}_1\theta,\\
 \text{3D}: \ &u=(-\Delta)^{-1} \left(-\mathcal{R}_1\mathcal{R}_3,-\mathcal{R}_2\mathcal{R}_3,\mathcal{R}_1^2+\mathcal{R}_2^2\right)\theta.
 \end{split}
\end{equation*}
We also define the Stokes-$\alpha$ model
\eqn \label{Stokesalpha}
\theta_{t}+u\cdot\nabla\theta=0, \quad u=\Lambda^{2\alpha-2}\mathcal{R}^{\perp}\mathcal{R}_1\theta
\een 
which  interpolates between \eqref{Stokes} ($\alpha=0$) and \eqref{IPM} ($\alpha=1$).

\subsubsection*{Magnetogeostrophic equation}
This equation is of the form
\eqn \label{MG}
\theta_{t}+u\cdot\nabla\theta=0, \quad \hat{u}(\xi)=\frac{\left(\xi_2\xi_3|\xi|^{2}-\xi_1\xi_2^2\xi_3,-\xi_1\xi_3|\xi|^{2}-\xi_2^3\xi_3,\xi_1^2\xi_2^2+\xi_2^4\right)}{|\xi|^{2}\xi_3^2+\xi_2^4}\hat{\theta}(\xi),
\een 
where $\theta$ represents the buoyancy field. This model is proposed by Moffatt and Loper  \cite{Moffatt} as a reduced model of the full magnetohydrodymics to study the geodynamo and turbulence in the Earth's fluid core \cite{Friedlander2, Friedlander3, Moffatt}.

\subsubsection*{Patch problems}
There is a large (and growing) literature in a particular class of weak solution with a \emph{patch} type initial data:
$$
\theta_0(x)=\theta^1\textbf{1}_{\Omega^1}+\theta^2\textbf{1}_{\mathbb{R}^n\setminus\Omega^1},
$$
with $\theta^i$ positive constants. In the case of the surface quasi-geostrophic equation, these initial data correspond to a sharp front between two different temperatures and they have obvious interest in meteorology. For the incompressible porous medium equation this situation is known as Muskat problem. This problem is of practical interest because it is a model of the dynamics of the interface between two different fluids in oil wells or geothermal reservoirs. We refer the readers to \cite{BCG, castro2012breakdown, ccfgl, ccgs-10, c-g07, CGO, gancedo2008existence, GG, G, RodrigoQG} and references therein for more details.

All the previous models are posed in two or three spatial dimensions. However, several related one dimensional problems have been studied. The one dimensional reduction idea was initiated by Constatin-Lax-Majda \cite{Constantin}; they proposed the following 1D model 
\[
\omega_{t}=\omega\mathcal{H}\omega
\]
for the 3D Euler equation in the vorticity form and proved that $\omega$ blows up in finite time under certain conditions. Motivated by this work, other similar models were proposed and analyzed in the literature \cite{Baker, Chae, Cordoba 2, De Gregorio, Morlet}.

In this paper, we study an one dimensional model of the surface quasi-geostrophic equation and incompressible porous media equation, and one dimensional and $n$ dimensional models of the dissipative quasi-geostrophic equations and the dissipative incompressible porous medium equation in the periodic domain.  We begin with the following one dimensional model
\eqn \label{1d qg}
\theta_{t} +\left(\theta \mathcal{H}\theta\right)_{x} =0, 
\een
which is obtained by replacing the Riesz transforms by the Hilbert transform $\mathcal{H}$ in the divergence form of (\ref{SQG}) and \eqref{IPM}. We note that (\ref{1d qg}) is also proposed as a model of dislocation dynamics \cite{dislocation1, dislocation2, dislocation3, dislocation4} where $\theta$ is related to the density of fractures per length in the material. Equation \eqref{1d qg} and related models have been studied by different authors. In a series of papers by A. Castro, D. Chae, A. C\'ordoba, D. C\'ordoba and M. Fontelos, the authors addressed the well-posedness and finite time singularities \cite{CC, CC2, Chae, Cordoba 2}. In particular, A. Castro and D. C\'ordoba proved global well-posedness for positive $L^2\cap C^{0,\gamma}$, vanishing at infinity data and finite time singularities for initial data such that $\theta_0(x_0)=0$ for some $x_0$. J. Carrillo, L. Ferreira, and J. Precioso in \cite{Carrillo} proved the existence of solution corresponding to initial data that are probability measures with finite second moment using gradient flows tools.

We next consider a dissipative model of (\ref{1d qg}):  
\eqn \label{1d dqg}
\theta_{t} +(1-\delta)\mathcal{H}\theta\theta_{x} +\delta (\theta \mathcal{H}\theta)_x + \nu \Lambda^{\gamma}\theta=0, \quad 0\leq\delta\leq1.
\een
We note that when $\nu=0$ and $\delta=1$, we return to (\ref{1d qg}). So, we can understand  \eqref{1d qg} as the limiting case of \eqref{1d dqg}. We note that there are several singularity formation results when $\nu=0$: $0<\delta<1/3$ and $\delta=1$ \cite{Morlet}, $0<\delta\leq1$  \cite{Chae}, and $\delta=0$ \cite{Cordoba 2, LiRodrigo2}. The case $\delta=-1$ is similar to the Kuramoto-Sivashinsky equation and it has been studied in \cite{LiRodrigo}.

Finally, we analyze two $n-$ dimensional versions of (\ref{1d dqg}). The first model is the equation with $\delta>0$ and $\nu=0$:
\eqn \label{nd dqg}
\theta_{t} +(1-\delta)u\cdot\nabla\theta + \delta \nabla\cdot(\theta \mathcal{R}\theta)=0, \quad \delta>0.
\een
The second model is the equation with $\delta=0$ and $\nu>0$ which corresponds to the dissipative quasi-geostrophic equation:
\eqn \label{nd dqg 2}
\theta_{t} +u\cdot\nabla\theta + \nu \Lambda^{\gamma}\theta=0.
\een
Here, $u$ in (\ref{nd dqg}) and (\ref{nd dqg 2}) is a divergence-free vector field $u$ such that
\[
\widehat{u}(k)=m(k) \widehat{\theta}(k), \quad k\cdot m(k)=0, \quad m\in L^{\infty}.
\]

Compared with existing results showing global solutions or blowups in finite time with smooth initial data, we establish several global existence of weak solutions with rough ($L^{1+s},s>0$) initial data. To the best of our knowledge, the sharpest global existence result is for initial data in $L^p$ requires $p>4/3$ \cite{Marchand}. To this end, we carefully choose dissipative quantities to minimize conditions of initial data. These quantities have the same flavor as the \emph{Shannon entropy} 
\[
\int \theta\log \theta dx.
\]
For more applications of the entropy, see \cite{blanchet2010functional, Dolbeault3, carrillo2012uniqueness}, where the authors apply these ideas jointly with the logarithmic Hardy-Littlewood-Sobolev inequality to the parabolic-elliptic Keller-Segel equation in the plane.

\section{Preliminaries}
In this paper, all constants will be denoted by $C$ that is a generic constant depending only on the quantities specified in the context.

\subsection*{Hilbert Transform}
The Hilbert transfom is defined as 
\[
\mathcal{H}f(x)=\frac{1}{2\pi} \text{p.v.} \int_{\mathbb{T}} \frac{f(y)}{\tan\left((x-y)/2\right)}dy.
\]
We list several properties of the Hilbert transform in $\mathbb{T}$. 
\begin{equation} \label{hilbert}
 \begin{split}
 & \mathcal{H} \left(\mathcal{H}f\right)=-f + \langle f \rangle, \quad \langle f \rangle=\int_{\mathbb{T}}fdx, \\
 & \left(\mathcal{H}f\right)_{x}=\mathcal{H}(f_{x}),\quad \langle \mathcal{H}f \rangle=0\\
 & \mathcal{H}\left(f\mathcal{H}g+ g\mathcal{H}f\right)=\left(\mathcal{H}f\right)\left(\mathcal{H}g\right)-fg - \langle f \rangle \langle g \rangle, \\
 & \int_{\mathbb{T}} \left(\mathcal{H}f\right) g dx=-\int_{\mathbb{T}} f\left(\mathcal{H}g\right) dx.
 \end{split}
\end{equation}

\subsection*{Function Spaces}
The energy norms in $\mathbb{T}^{n}$ are defined as follows: 
\begin{equation} \label{energy}
 \begin{split}
 & \|f\|^{2}_{L^{2}(\mathbb{T}^{n})}=\sum_{k\in \mathbb{Z}^{n}}  \left|\widehat{f}(k)\right|^{2}, \quad  \|f\|^{2}_{\dot{H}^{s}(\mathbb{T}^{n})}=\sum_{k\in \mathbb{Z}^{n}\setminus\{0\}} |k|^{2s}\left|\widehat{f}(k)\right|^{2} \\
 &  \|f\|^{2}_{H^{s}(\mathbb{T}^{n})}=\sum_{k\in \mathbb{Z}^{n}} \left(1+|k|^{s}\right)^{2} \left|\widehat{f}(k)\right|^{2}.
 \end{split}
\end{equation}
We also introduce the semi-norm
\[
\left\|f_{x}\right\|_{l^{1}}:=\sum_{k\in \mathbb{Z}}|k| \left| \widehat{f}(k)\right|,
\]
that is, the derivative is in the Wiener algebra of absolutely convergent Fourier series.

\subsection*{Operator $\Lambda^{\gamma}$}
The differential operator $\Lambda^{\gamma}=(\sqrt{-\Delta})^{\gamma}$ is defined by the action of the following kernels \cite{Cordoba}:
\eqn \label{lambda gamma}
\Lambda^{\gamma} f(x)=c_{\gamma,n}\text{p.v.} \int_{\mathbb{T}^n} \frac{f(x)-f(y)}{|x-y|^{n+\gamma}}dy+c_{\gamma,n}\sum_{k\in\mathbb{Z}^n\setminus\{0\}} \int_{\mathbb{T}^n} \frac{f(x)-f(y)}{|x-y+2k\pi|^{n+\gamma}}dy,
\een
where $c_{\gamma,n}>0$ is a normalized constant. In particular, in one dimension with $\gamma=1$,  
\[
\Lambda f(x)=\mathcal{H}f_{x}(x)=\frac{1}{2\pi} \text{p.v.} \int_{\mathbb{T}} \frac{f(x)-f(y)}{\sin^{2}\left((x-y)/2\right)}dy
\]

\subsection*{Minimum Principle}
In this paper, we assume that $\theta_{0}\ge 0$. It is well-known that this property propagates for (\ref{1d qg}), (\ref{1d dqg}), (\ref{nd dqg}) and  (\ref{nd dqg 2}). Let us give a sketch of the proof for (\ref{nd dqg 2}) \cite{Chae, Cordoba}. (The same argument holds for (\ref{1d qg}), (\ref{1d dqg}), (\ref{nd dqg}).) Let's assume that $\theta(x,t)\in C^1([0,T]\times\mathbb{T}^n)$ and $x_t$ be a point such that $m(t)=\min_x \theta(x,t)=\theta(x_t,t)$. Since $m(t)$ is a continuous Lipschitz function, it is differentiable at almost every point $t$ by Rademacher's theorem. Then, from the definition of $\Lambda^{\gamma}$ and the non-negative assumption, we have 
\[
\Lambda^{\gamma} \theta(x_t,t)= c_{\gamma,n}\text{p.v.} \int_{\mathbb{T}^n} \frac{\theta(x_t)-\theta(x_t-y)}{|y|^{n+\gamma}}dy+c_{\gamma,n}\sum_{k\in\mathbb{Z}^n\setminus\{0\}} \int_{\mathbb{T}^n} \frac{\theta(x_t)-\theta(x_t-y)}{|y-2k\pi|^{n+\gamma}}dy\leq 0.
\]
This implies that 
\[
m^{'}(t)=-\nu \Lambda^{\gamma} \theta(x_t,t) \ge 0. 
\]
Therefore, we conclude that $\theta(t,x)\ge 0$ for all time. In the paper we will deal with weak solution that are not continuous in general. However, for the regularized problems, $\theta^\epsilon$, the same argument works. Then we construct $\theta$ as the limit (in the appropriate space) of $\theta^\epsilon$. As the $\theta$ will be also the pointwise limit of $\theta^\epsilon$, we conclude.

\subsection*{Compactness} 
Since we look for weak solutions, we use following compactness arguments when we pass to the limit in weak formulations.

\begin{lemma}[\cite{Temam}] \label{lemma:2.1}
Let $X_0,X,X_1$ be reflexive Banach spaces such that
\[
X_0\subset  X\subset X_1,
\] 
where $X_0$ is compactly embedded in $X$. Let $T>0$ be a finite number and let $\alpha_0$ and $\alpha_1$ be two finite numbers such that $\alpha_i>1$. Then, 
$$
Y=\left\{u\in L^{\alpha_0}\left(0,T; X_0\right),\  \partial_t u\in L^{\alpha_1}\left(0,T; X_1\right)\right\}
$$
is compactly embedded in $L^{\alpha_0}\left(0,T; X\right)$.
\end{lemma}

\begin{lemma} [\cite{Lions1}] \label{lemma:2.2}
Let $\Omega$ be a bounded set in $\mathbb{R}^{n}$. Let $\left(g^{\epsilon}\right)$ and $\left(h^{\epsilon}\right)$ converge weakly to $g$ and $h$ respectively in $L^{p_{1}}\left(0,T; L^{p_{2}}(\Omega)\right)$ and $L^{q_{1}}\left(0,T; L^{q_{2}}(\Omega)\right)$, with
\[
1 \leq p_{1}, p_{2} \leq \infty, \quad \frac{1}{p_{1}}+\frac{1}{q_{1}}=\frac{1}{p_{2}}+\frac{1}{q_{2}}=1.
\]
Suppose that we have the following properties uniformly in $\epsilon>0$:
\begin{equation} \label{compact conditions} 
 \begin{split}
 & g^{\epsilon}_{t} \ \text{is bounded in } \  L^{1}\left(0,T; W^{-m,1}(\Omega)\right) \  \text{for some} \ m\ge 0,\\
 & \left\|h^{\epsilon}-h^{\epsilon}(\cdot +y,t)\right\|_{L^{q_{1}}(0,T; L^{q_{2}}(\Omega))} \rightarrow 0 \ \text{as} \ |y| \rightarrow 0.
  \end{split}
\end{equation}
Then, $\left(g^{\epsilon}h^{\epsilon}\right)$ converges to $gh$ in the sense of distributions.
\end{lemma}

\noindent
We note that the second condition in (\ref{compact conditions}) holds when $\left(h^{\epsilon}\right)$ has positive regularity in space.

\section{Statements of Results}

\subsection{1D model of \eqref{SQG} and \eqref{IPM}}
We consider the equation \eqref{1d qg} in $\mathbb{T}$ with non-negative initial data. Since (\ref{1d qg}) satisfies the minimum principle, $\theta(t,x)\ge 0$ for all time. We notice that (\ref{1d qg}) is dissipative: the entropy
\eqn \label{entropy}
\mathfrak{E}(\theta)=\int_{\mathbb{T}} \left[\theta \log\theta -\theta +1\right] dx
\een
gains $\dot{H}^{1/2}$ regularity. Therefore, it is natural to assume that $\theta_{0}$ satisfies the following conditions
\eqn \label{1d qg initial data}
\theta_{0}(x)\ge 0, \quad  \theta_{0}(x)\in L^{1+s},s>0, \quad \mathfrak{E}(\theta_{0})<\infty,
\een
where we need the second condition to obtain a uniform bound of $\mathfrak{E}(\theta^{\epsilon}_{0})$ when we construct an approximate sequence of solutions. We define the function space as follows:
\begin{equation}\label{Aspace} 
 \begin{split}
 \mathcal{A}_T&=\left\{\theta \in L^{\infty}\left(0,T; L^{1}(\mathbb{T})\right): \sup_{0\leq t<T} \mathfrak{E}(\theta(t)) + \int^{T}_{0}\left\|\Lambda^{1/2}\theta(t)\right\|^{2}_{L^{2}}dt<\infty \right\}.
  \end{split}
\end{equation}

\begin{definition} 
$\theta$ is a weak solution of (\ref{1d qg}) if $\theta\in \mathcal{A}_T$ and (\ref{1d qg}) holds in the sense of distributions: for any $\psi\in \mathcal{C}^{\infty}_{c}\left([0,T)\times\mathbb{T}\right)$, 
\[
\int^{T}_{0}\int_{\mathbb{T}} \Big[\theta \psi_{t} + \theta \left(\mathcal{H}\theta\right) \psi_{x}\Big]  dxdt =\int_{\mathbb{T}}\theta_{0}(x)\psi(x,0)dx \quad \text{for any $T<\infty$}.
\]
\end{definition}

\begin{theorem}\label{thm:3.1} 
For any initial datum $\theta_{0}$ satisfying (\ref{1d qg initial data}), there exists a weak solution of (\ref{1d qg}) in $\mathcal{A}_T$ for all $T>0$. 
\end{theorem}

\subsection{Dissipative 1D model}
We first show the local well-posedness of \eqref{1d dqg} in $H^{2}(\mathbb{T})$ with $\nu>0$ and $\delta=0$. Then, we continue to prove that the solution can be extended beyond $T>0$ as long as $\left\|\theta_{x}(t)\right\|_{L^{\infty}}$ is integrable in $[0,T]$.  The last condition can be achieved if $\left\|\theta_{0x}\right\|_{l^{1}}$ is sufficiently small. We note that H. Dong and A. Kiselev proved the global existence in the critical case $\gamma=1$ in \cite{HDong} and \cite{Kiselev}, respectively. Therefore, we restrict ourselves to the case $0< \gamma<1$ for the local well-posedness.

\begin{theorem}\label{thm:3.2}
Let $\nu> 0$, $\delta= 0$, and $0< \gamma<1$. For any initial datum $\theta_{0} \in H^{2}(\mathbb{T})$, there exists $T=T(\theta_{0})>0$ such that there exists a unique solution of (\ref{1d dqg}) $\theta \in C\left(0,T; H^{2}(\mathbb{T})\right)$. If  $\displaystyle \left\|\theta_{0x}\right\|_{l^{1}}<\nu$, we can take $T=\infty$. 
\end{theorem}

We recall that the local existence part stated in this theorem was proved also in \cite{HDong}.

The second result is the global existence of a weak solution of (\ref{1d dqg}). We note that the additional term in the right-hand side of (\ref{1d dqg}) depletes the nonlinear term $\left(\mathcal{H}\theta\right) \theta_{x}$ when $\delta\ge \frac{1}{2}$ with a strictly positive lower bound of $\theta_{0}$. This enables us to show the existence of a weak solution with $\nu \ge 0$. However, regularity of initial data prescribed below is relatively higher than the usual $L^{2}$ regularity in dissipative equations because (\ref{1d dqg}) is not in the divergence form.  In this paper, we assume that initial data satisfy the following conditions:
\eqn \label{1d dqg initial data}
\theta_{0}\in H^{1/2}(\mathbb{T}), \quad m_{0}=\text{ess\,inf}_{x\in \mathbb{T}}\theta_{0}(x)>0.
\een

\begin{definition}
$\theta$ is a weak solution of (\ref{1d dqg}) with $\theta_{0} \in H^{1/2}(\mathbb{T})$ if \ $\theta(t) \in H^{1/2}(\mathbb{T})$ for any $t\leq T$ and (\ref{1d dqg}) holds in the sense of distributions: for any $\psi\in \mathcal{C}^{\infty}_{c}\left([0,T)\times \mathbb{T}\right)$, 
\[
\int^{T}_{0}\int_{\mathbb{T}} \big[\theta \psi_{t} -(1-\delta)\mathcal{H}\theta\theta_{x} \psi -\nu\theta \Lambda^{\gamma}\psi+\delta \theta \mathcal{H}\theta\psi_{x} \big] dxdt= \int_{\mathbb{T}}\theta_{0}(x)\psi(x,0)dx \quad \text{for any $T<\infty$}.
\]
\end{definition}

\begin{theorem}\label{thm:3.3}
Let $\nu\ge 0$, $\gamma\ge 0$ and $1/2 \leq \delta<1$. For any  initial datum $\theta_{0}$ satisfying (\ref{1d dqg initial data}), there exists a global weak solution of \eqref{1d dqg} such that 
\[
\theta \in L^{\infty}\left(0,\infty; H^{1/2}(\mathbb{T})\right) \cap L^{2}\left(0,\infty; H^{\max\{1,(1+\gamma)/2\}}(\mathbb{T})\right).
\]
Moreover, such a solution is unique in $L^{2}(\mathbb{T})$ if $\nu>0$ and $\gamma\ge 2$. 
\end{theorem}

\subsection{High dimensional model}
We finally consider the equation \eqref{nd dqg} and \eqref{nd dqg 2} in $\mathbb{T}^n$, $n=2,3$, with a divergence-free vector field $u$ satisfying 
\eqn \label{hypu}
\widehat{u}(k)=m(k) \widehat{\theta}(k), \quad k\cdot m(k)=0, \quad m\in L^{\infty}.
\een
We begin with the equation \eqref{nd dqg}. We use the $n$-dimensional version of the  entropy \eqref{entropy} and functional space \eqref{Aspace}. This is due to the fact that the advection term vanishes in the computation of $\mathfrak{E}(\theta)_{t}$ by the divergence-free condition of $u$.

\begin{definition} 
$\theta$ is a weak solution of \eqref{nd dqg} if $\theta\in \mathcal{A}_T$ and \eqref{nd dqg} holds in the sense of distributions: for any $\psi\in \mathcal{C}^{\infty}_{c}\left([0,T)\times\mathbb{T}^n\right)$
\[
\int^{T}_{0}\int_{\mathbb{T}^n} \Big[\theta \psi_{t} + \theta u\cdot\nabla\psi+\theta \mathcal{R}\theta\cdot\nabla\psi\Big]  dxdt =\int_{\mathbb{T}^n}\theta(x,0)\psi(x,0)dx \quad \text{for every $T<\infty$.}
\]
\end{definition}

\begin{theorem}\label{thm:3.4}
For any initial datum $\theta_{0}$ satisfying \eqref{1d qg initial data}, there exists a  weak solution of \eqref{nd dqg}  in $\mathcal{A}_T$ for all $T>0$.  
\end{theorem}

We note that the same result holds for a smoother velocity field: 
\begin{equation}\label{hypu2}
\widehat{u}(k)=|k|^\beta m(k) \widehat{\theta}(k), \quad k\cdot m(k)=0, \quad m(0)=0
\end{equation}
with bounded $m$ and $\beta<0$. We note that $u$ in (\ref{hypu2}) covers  \eqref{Euleralpha},  \eqref{Stokes} and \eqref{Stokesalpha}.

\begin{corollary}\label{cor:3.5}
For any initial datum $\theta_{0}$ satisfying \eqref{1d qg initial data} and $u$ given by \eqref{hypu2}, there exists a weak solution of \eqref{nd dqg} in $\mathcal{A}_T$ for every $T>0$. 
\end{corollary}

We finally deal with  the equation \eqref{nd dqg 2} with a slightly different entropy
\[
\mathfrak{E}(\theta)=\int_{\mathbb{T}^{n}} (\theta+1) \log(\theta+1) dx.
\]
Due to the lack of the smoothing effect from $\nabla \cdot (\theta \mathcal{R}\theta)$, we only have bounds of $\theta$ for $t\ge \tau>0$. Therefore, we define the function space and the notion of weak solution as follows.
\begin{equation}\label{Bspace} 
 \begin{split}
 \mathcal{B}_T&=\bigg\{\theta \in L^{\infty}\left(0,T; L^{1}\left(\mathbb{T}^{n}\right)\right) \cap L^{\infty}\left([\tau,T); L^{\infty}\left(\mathbb{T}^{n}\right)\right):\\
 & \sup_{\tau\leq  t<T} \mathfrak{E}(\theta(t)) + \int^{T}_{\tau}\left\|\Lambda^{\gamma/2}\theta(t)\right\|^{2}_{L^{2}}dt<\infty \bigg\}.
  \end{split}
\end{equation}

\begin{definition} 
$\theta$ is a weak solution of \eqref{nd dqg 2} if $\theta\in \mathcal{B}_T$ and \eqref{nd dqg} holds in the sense of distributions: for any $\psi\in \mathcal{C}^{\infty}_{c}\left([\tau,T)\times\mathbb{T}^n\right)$
\[
\int^{T}_{\tau}\int_{\mathbb{T}^n} \Big[\theta \psi_{t} + \theta u\cdot\nabla\psi+\nu\theta \Lambda^{\gamma}\psi\Big]  dxdt =0 \quad \text{for every $0<\tau<T<\infty$.}
\]
\end{definition}

\begin{theorem}\label{thm:3.6}
For any initial datum $\theta_{0}$ satisfying \eqref{1d qg initial data}, there exists a  weak solution of \eqref{nd dqg 2}  in $\mathcal{B}_T$ for every $T>0$.  Moreover, $\theta(t)$ converges to $\theta_{0}$ in $H^{-2}(\mathbb{T}^{n})$ as $t\rightarrow 0$.
\end{theorem}

\begin{remark} 
Actually, following the ideas in the proof of Theorem \ref{thm:3.6}, we can prove that the solution $\theta$ in Theorems \ref{thm:3.1}-\ref{thm:3.6} is in $L^\infty(\tau,T;L^\infty(\mathbb{T}^n))$ for every $0<\tau<T<\infty$.
\end{remark}

The proofs of our results are outlined as follows. We first obtain a priori estimates in given function spaces. We then generate approximate sequence of solutions and pass to the limits in weak formulation using Lemma \ref{lemma:2.1} or \ref{lemma:2.2}.

\section{Proof of Theorem \ref{thm:3.1}} \label{sec:4}

\noindent 
We consider the equation (\ref{1d qg}) 
\[
\theta_{t} +\left(\theta \mathcal{H}\theta\right)_{x} =0
\]
and the entropy
\[
\mathfrak{E}(\theta)=\int_{\mathbb{T}}\left[\theta \log\theta-\theta +1\right]dx.
\]
Since $\theta(t)\ge 0$, $\mathfrak{E}(\theta)\ge 0$. Moreover, the direct computation yields that 
\begin{equation} \label{eq:4.1}
 \begin{split}
 \frac{d}{dt}\mathfrak{E}(\theta)& =\int_{\mathbb{T}} \left[\theta_{t} \log \theta(t)+\theta(t)(\log \theta(t))_{t} -\theta_{t}\right]dx =\int_{\mathbb{T}}\theta_{t} \log \theta(t)dx \\
 & =-\int_{\mathbb{T}} \left(\theta \mathcal{H}\mathcal{\theta}\right)_{x} \log \theta dx =\int_{\mathbb{T}}\left(\mathcal{H}\theta\right)\theta_{x}dx=-\int_{\mathbb{T}}\theta \Lambda \theta dx=-\left\|\Lambda^{1/2}\theta\right\|^{2}_{L^{2}}.
  \end{split}
\end{equation}
Therefore, we have $\theta \in \mathcal{A}$. We now construct a sequence of solutions $\left(\theta^{\epsilon}\right)$ by solving
\[
\theta^{\epsilon}_{t} +\left(\theta^{\epsilon}\mathcal{H}\theta^{\epsilon}\right)_{x} =\epsilon \theta^{\epsilon}_{xx},\quad  \theta^{\epsilon}_{0}=\rho_{\epsilon}\ast \theta_{0},
\]
where $\rho_{\epsilon}$ is a standard mollifier. Then, $\theta^{\epsilon}$ satisfies that
\eqn \label{eq:4.2}
\frac{d}{dt}\mathfrak{E}(\theta^{\epsilon}) +\left\|\Lambda^{1/2}\theta^{\epsilon}\right\|^{2}_{L^{2}} +4\epsilon\int_{\mathbb{T}}\left|\left(\sqrt{\theta^{\epsilon}}\right)_{x}\right|^{2}dx=0.
\een
Integrating (\ref{eq:4.2}) in time, we have
\begin{equation} \label{eq:4.3}
 \begin{split}
 \mathfrak{E}(\theta^{\epsilon}(t))  +\int^{t}_{0}\left\|\Lambda^{1/2}\theta^{\epsilon}(s)\right\|^{2}_{L^{2}}ds  + 4\epsilon \int^{t}_{0}\int_{\mathbb{T}}\left|\left(\sqrt{\theta^{\epsilon}}(s)\right)_{x}\right|^{2}dxds = \mathfrak{E}(\theta^{\epsilon}_{0}).
  \end{split}
\end{equation}
Since $x\log x -x+1 \leq  x^{s+1}+1$ for $x\ge 0$, we can bound the last term in (\ref{eq:4.3}) as 
\[
\mathfrak{E}(\theta^{\epsilon}_{0})  \leq 2\pi + \left\| \theta^{\epsilon}_{0}\right\|^{s+1}_{L^{s+1}} \leq 2\pi + \left\| \theta_{0}\right\|^{s+1}_{L^{s+1}}.
\]
Therefore, the sequence $\left(\theta^{\epsilon}\right)$ is uniformly bounded in $\mathcal{A}_T$. By Poincar\'e's inequality, we obtain uniform bounds of $\theta^{\epsilon}$  and $\mathcal{H}\theta^{\epsilon}$ in $L^{2} \left((0,T); H^{1/2}(\mathbb{T})\right)$.  Moreover, by interpolating $L^{\infty} \left(0,T; L^{1}(\mathbb{T})\right)$ and $L^{2} \left((0,T); H^{1/2}(\mathbb{T})\right)$, we  have uniform bounds of $\theta^\epsilon$ and $\mathcal{H}\theta^\epsilon$ in $L^4\left(0,T; L^2(\mathbb{T})\right)$. These estimates and the duality pairing imply that 
\begin{equation*} 
 \begin{split}
 \theta^{\epsilon}_{t} =-\left(\theta^{\epsilon}\mathcal{H}\theta^{\epsilon}\right)_{x} +\epsilon \theta^{\epsilon}_{xx} \in L^{2}\left(0,T; H^{-2}(\mathbb{T})\right)
  \end{split}
\end{equation*}
uniformly in $\epsilon>0$. Lemma \ref{lemma:2.1} with
\[
X_0=L^2\left(0,T; H^{1/2}(\mathbb{T})\right),\quad X=L^2\left(0,T; L^2(\mathbb{T})\right),\quad X_1=L^2\left(0,T; H^{-2}(\mathbb{T})\right),
\]
allows to pass to the limit in
\[
\int^{T}_{0}\int_{\mathbb{T}} \Big[\theta^{\epsilon} \psi_{t} + \theta^{\epsilon} \left(\mathcal{H}\theta^{\epsilon}\right) \psi_{x}\Big]  dxdt =\int_{\mathbb{T}}\theta^\epsilon_{0}(x)\psi(x,0)dx
\]
to obtain a weak solution in $\mathcal{A}_T$.

\section{Proof of Theorem \ref{thm:3.2}}\label{sec:5}

\subsection{Local Well-posedness} 
We consider the equation (\ref{1d dqg}) with $\delta=0$: 
\eqn \label{eq:5.1}
\theta_{t} +\left(\mathcal{H}\theta\right)\theta_{x} + \nu \Lambda^{\gamma}\theta=0, \quad \nu>0, \ \ \gamma>0. 
\een
We here only provide a priori estimates. We first multiply (\ref{eq:5.1}) by $\theta$ and integrate over $\mathbb{T}$:
\begin{equation*} 
 \begin{split}
 \frac{1}{2}\frac{d}{dt}\left\| \theta \right\|^{2}_{L^{2}}+\nu\left\|\Lambda^{\frac{\gamma}{2}}\theta\right\|^{2}_{L^{2}}&=-\int_{\mathbb{T}}\Big[\left(\mathcal{H}\theta\right)\theta_{x}\theta\Big]dx\leq \|\theta_{x}\|_{L^{\infty}} \left\| \mathcal{H}\theta\right\|_{L^{2}} \|\theta\|_{L^{2}} \leq \|\theta_{x}\|_{L^{\infty}} \|\theta\|^{2}_{L^{2}}. 
 \end{split}
\end{equation*}
We next take $\partial_{x}$ to (\ref{eq:5.1}) and do the energy estimate. 
\begin{equation*} 
 \begin{split}
 &\frac{1}{2}\frac{d}{dt}\left\|\theta_{x} \right\|^{2}_{L^{2}}+\nu\left\|\Lambda^{\frac{\gamma}{2}}\theta_{x}\right\|^{2}_{L^{2}} =-\int_{\mathbb{T}} \left(\mathcal{H}\theta \theta_{x}\right)_{x}\theta_{x}dx=-\int_{\mathbb{T}}\left(\mathcal{H}\theta\right)_{x}\left(\theta_{x}\right)^{2}dx-\int_{\mathbb{T}} \mathcal{H}\theta \theta_{xx}\theta_{x}dx \\
 &=-\int_{\mathbb{T}}\left(\mathcal{H}\theta\right)_{x}\left(\theta_{x}\right)^{2}dx-\frac{1}{2}\int_{\mathbb{T}} \mathcal{H}\theta \left[\left(\theta_{x}\right)^{2}\right]_{x}dx =-\frac{1}{2}\int_{\mathbb{T}}\left(\mathcal{H}\theta\right)_{x}\left(\theta_{x}\right)^{2}dx \\
 &\leq \|\theta_{x}\|_{L^{\infty}}  \left\| \mathcal{H}\theta_{x}\right\|_{L^{2}}  \|\theta_{x}\|_{L^{2}} \leq \|\theta_{x}\|_{L^{\infty}}\|\theta_{x}\|^{2}. 
 \end{split}
\end{equation*}
Similarly, by taking $\partial_{xx}$ to (\ref{eq:5.1}), we have 
\[
\frac{1}{2}\frac{d}{dt}\left\|\theta_{xx} \right\|^{2}_{L^{2}}+\nu\left\|\Lambda^{\frac{\gamma}{2}}\theta_{xx}\right\|^{2}_{L^{2}}\leq  \left(\left\|\theta_{x}\right\|_{L^{\infty}} + \left\|\mathcal{H}\theta_{x}\right\|_{L^{\infty}} \right) \|\theta_{xx}\|^{2}_{L^{2}}. 
\]
Therefore, we obtain that 
\begin{equation} \label{eq:5.2}
 \begin{split}
 \frac{d}{dt} \left\| \theta \right\|^{2}_{H^{2}} +\nu\left\|\Lambda^{\frac{\gamma}{2}}\theta\right\|^{2}_{H^{2}} \leq C\left(\left\|\theta_{x}\right\|_{L^{\infty}} + \left\|\mathcal{H}\theta_{x}\right\|_{L^{\infty}} \right) \|\theta\|^{2}_{H^{2}}. 
 \end{split}
\end{equation}
Since $\left\|\theta_{x}\right\|_{L^{\infty}} + \left\|\mathcal{H}\theta_{x}\right\|_{L^{\infty}} \leq C\|\theta\|_{H^{2}}$, (\ref{eq:5.2}) implies the local well-posedenss in $H^{2}(\mathbb{T})$. Moreover, by integrating (\ref{eq:5.2}) in time, we obtain that 
\[
\left\| \theta(t)\right\|^{2}_{H^{2}} \leq C\left\| \theta_{0}\right\|^{2}_{H^{2}} \exp \int^{t}_{0} C\big(\|\theta_{x}(s)\|_{L^{\infty}} + \left\|\mathcal{H}\theta_{x}(s)\right\|_{L^{\infty}}\big) ds.
\]
Using the logarithmic bound (similar to the Beale-Kato-Majda criterion \cite{Beale})
\[
\|\mathcal{H} \theta_x\|_{L^\infty}\leq C(1+\|\theta_x\|_{L^\infty}\log(e+\|\theta\|_{H^2})+\|\theta_x\|_{L^2}),
\]
the solution can be continued as long as we can control $\left\|\theta_{x}\right\|_{L^{\infty}}$.

\subsection{Estimation of $\left\|\theta_{x}\right\|_{l^{1}}$} 
We now control $\left\|\theta_{x}\right\|_{L^{\infty}}$ by $\left\|\theta_{x}\right\|_{l^{1}}$. For $k \in \mathbb{Z}$,
\begin{equation*} 
 \begin{split}
 |k| \left|\widehat{\theta}(k)\right|_{t} &=- \nu|k|^{1+\gamma} \left|\widehat{\theta}(k)\right| - \frac{\widehat{\theta}(k)}{\left|\widehat{\theta}(k)\right|} |k| \sum_{l\in \mathbb{Z}}\left[i l\widehat{\theta}(l)\frac{i(k-l)}{|k-l|}\widehat{\theta}(k-l)\right] \\
 &\leq - \nu|k|^{1+\gamma} \left|\widehat{\theta}(k)\right| + |k| \sum_{l\in \mathbb{Z}} \left|\widehat{\theta}(k-l) \right| |l| \left|\widehat{\theta}(l)\right|.
 \end{split}
\end{equation*}
By taking the summation over $k\in \mathbb{Z}$, 
\begin{equation*}
 \begin{split}
 \frac{d}{dt} \left\| \theta_{x} \right\|_{l^{1}} &\leq -\nu\sum_{k\in \mathbb{Z}} |k|^{1+\gamma} \left| \widehat{\theta}(k)\right|  +  \sum_{k\in \mathbb{Z}} \sum_{l\in \mathbb{Z}}  |k| |l| \left|\widehat{\theta}(l) \right| \left|\widehat{\theta}(k-l)\right| \\
 & = -\nu\sum_{k\in \mathbb{Z}} |k|^{1+\gamma} \left|\widehat{\theta}(k)\right|  + \sum_{l\in \mathbb{Z}}  |l| \left|\widehat{\theta}(l) \right| \sum_{k\in \mathbb{Z}} |k| \left|\widehat{\theta}(k-l)\right| \\
 & \leq -\nu\left\|\Lambda^{1+\gamma}\theta \right\|_{l^{1}}+ \left\| \theta_{x} \right\|^{2}_{l^{1}}\leq \left\| \theta_{x} \right\|_{l^{1}} \left(\left\| \theta_{x} \right\|_{l^{1}} -\nu\right).
 \end{split}
\end{equation*}
Therefore, $\left\| \theta_{x}(t) \right\|_{l^{1}}<\nu$ as long as $\left\| \theta_{0,x} \right\|_{l^{1}}<\nu$. This completes the proof.

\section{Proof of Theorem \ref{thm:3.3}} \label{sec:6}
\noindent 
We consider the equation (\ref{1d dqg}) which is equivalent to 
\eqn \label{eq:6.1}
\theta_{t} +\left(\mathcal{H}\theta\right)\theta_{x} + \nu\Lambda^{\gamma}\theta+\delta \theta \Lambda \theta=0, \quad \nu\ge0, \  \gamma\ge 0, \ 1/2 \leq \delta<1.
\een
We begin with  a priori estimates. To obtain the $L^{2}$ bound, we multiply (\ref{eq:6.1}) by $\theta$ and integrate over $\mathbb{T}$:
\begin{equation*} 
 \begin{split}
 &\frac{1}{2}\frac{d}{dt}\left\| \theta \right\|^{2}_{L^{2}}+\nu\left\|\Lambda^{\gamma/2}\theta\right\|^{2}_{L^{2}} =-\int_{\mathbb{T}} \left[\left(\mathcal{H}\theta\right)\theta_{x}\theta\right] dx - \delta\int_{\mathbb{T}}\left[\theta^{2} \Lambda \theta\right]dx \\
 &=-\frac{1}{2}\int_{\mathbb{T}}\left[\left(\mathcal{H}\theta\right)\left(\theta^{2}\right)_{x}\right]dx - \delta\int_{\mathbb{T}}\left[\theta^{2} \Lambda \theta\right]dx =\left(\frac{1}{2}-\delta\right)\int_{\mathbb{T}}\left[\theta^{2} \Lambda \theta\right]dx.
  \end{split}
\end{equation*}
Since $\theta \ge 0$, we have 
\[
\int_{\mathbb{T}}\left[\theta^{2} \Lambda \theta\right]dx=\int_{\mathbb{T}} \int_{\mathbb{T}} \frac{\left(\theta(x)-\theta(y)\right)^{2}}{\sin^{2}\left((x-y)/2\right)} \cdot \frac{\theta(x)+\theta(y)}{2} dxdy\ge 0.
\]
Therefore, we obtain that
\eqn \label{eq:6.2}
\frac{1}{2}\frac{d}{dt}\left\| \theta \right\|^{2}_{L^{2}}+\nu\left\|\Lambda^{\gamma/2}\theta\right\|^{2}_{L^{2}}\leq 0. 
\een
We next obtain the $\dot{H}^{1/2}$ bound. We multiply (\ref{eq:6.1}) by $\Lambda \theta$ and integrate over $\mathbb{T}$:
\eqn \label{eq:6.3}
\frac{1}{2}\frac{d}{dt}\left\| \Lambda^{1/2}\theta \right\|^{2}_{L^{2}}+\nu\left\|\Lambda^{(1+\gamma)/2}\theta\right\|^{2}_{L^{2}}=-\int_{\mathbb{T}}\Big[\left(\mathcal{H}\theta\right)\theta_{x}\Lambda \theta\Big]dx - \delta\int_{\mathbb{T}}\Big[\theta \left(\Lambda \theta\right)^{2}\Big]dx.
\een
We now compute the first integral in the right-hand side of (\ref{eq:6.3}). Since  
\begin{equation*} 
 \begin{split}
 -\int_{\mathbb{T}}\left[\left(\mathcal{H}\theta\right)\theta_{x}\Lambda \theta\right]dx = \int_{\mathbb{T}}\left[\theta \mathcal{H} \left(\theta_{x}\left(\mathcal{H}\theta_{x}\right)\right)\right]dx =\frac{1}{2}\int_{\mathbb{T}}\left[\theta\left(\left(\Lambda \theta\right)^{2}-\left(\theta_{x}\right)^{2}\right)\right]dx,
  \end{split}
\end{equation*}
we have
\[
\frac{1}{2}\frac{d}{dt}\left\| \Lambda^{1/2}\theta \right\|^{2}_{L^{2}}+\nu\left\|\Lambda^{(1+\gamma)/2}\theta\right\|^{2}_{L^{2}}=\left(\frac{1}{2}-\delta\right)\int_{\mathbb{T}}\left[\theta\left(\Lambda \theta \right)^{2}\right]dx- \frac{1}{2}\int_{\mathbb{T}}\left[\theta \left(\theta_{x}\right)^{2}\right]dx,
\]
which implies that  
\eqn \label{eq:6.4}
\frac{1}{2}\frac{d}{dt}\left\| \Lambda^{1/2}\theta \right\|^{2}_{L^{2}}+\nu\left\|\Lambda^{(1+\gamma)/2}\theta\right\|^{2}_{L^{2}}+ \frac{1}{2}\int_{\mathbb{T}}\Big[\theta \left(\theta_{x}\right)^{2}\Big]dx\leq 0.
\een
We note that the minimum principle, with $m_{0}>0$, implies that $\theta(t,x)\ge m_{0}>0$ for all time. Therefore, by (\ref{eq:6.2}), (\ref{eq:6.4}), we obtain that 
\eqn \label{eq:6.5}
\frac{1}{2}\frac{d}{dt}\left\| \theta \right\|^{2}_{H^{1/2}}+\nu\left\|\Lambda^{\gamma/2}\theta\right\|^{2}_{H^{1/2}} + \frac{1}{2} m_{0} \left\| \theta_{x} \right\|^{2}_{L^{2}} \leq 0.
\een
Integrating (\ref{eq:6.5}) in time, we conclude that 
\begin{equation} \label{eq:6.6}
 \begin{split}
 & \theta \in L^{\infty}\left(0,\infty; H^{1/2}(\mathbb{T})\right) \cap L^{2}\left(0,\infty; \dot{H}^{1}(\mathbb{T})\right), \quad  \Lambda^{\gamma/2}\theta \in L^{2}\left(0,\infty; H^{1/2}(\mathbb{T})\right).
\end{split}
\end{equation}
We now construct an approximate sequence of solutions $(\theta^{\epsilon})$ by solving 
\[
\theta_{t} +\left(\mathcal{H}\theta\right)\theta_{x} + \Lambda^{\gamma}\theta= -\delta \theta \Lambda \theta+\epsilon \theta_{xx},\quad  \theta^{\epsilon}_{0}=\rho_{\epsilon}\ast \theta_{0}.
\]
Then, $(\theta^{\epsilon})$ is uniformly bounded in the space stated in (\ref{eq:6.6}). By interpolating $L^{\infty}\left(0,T; H^{1/2}(\mathbb{T})\right)$ and $L^{2}\left(0,T; H^{1}(\mathbb{T})\right)$, we have uniform bounds of $\mathcal{H}\theta^\epsilon$ and $\theta^{\epsilon}$ in $L^{4}\left(0,T; H^{3/4}(\mathbb{T})\right)$ that is embedded in $L^{4}\left(0,T; L^{2}(\mathbb{T})\right)$. This implies that 
\[
\left(\mathcal{H}\theta^{\epsilon}\right)\theta^{\epsilon}_{x} \in L^{4/3}\left(0,T; L^{1}(\mathbb{T})\right), \quad \theta^{\epsilon} \Lambda \theta^{\epsilon} \in L^{4/3}\left(0,T; L^{1}(\mathbb{T})\right)
\]
uniformly in $\epsilon>0$. Therefore, we obtain that 
\begin{equation*} 
 \begin{split}
 \theta^{\epsilon}_{xt} =-\Big(\left(\mathcal{H}\theta^{\epsilon}\right)\theta^{\epsilon}_{x} - \Lambda^{\gamma}\theta^{\epsilon}+\delta \theta^{\epsilon} \Lambda \theta^{\epsilon}+\epsilon \theta^{\epsilon}_{xx} \Big)_{x}\in L^{4/3}\left(0,T; W^{-2,1}(\mathbb{T})\right).
  \end{split}
\end{equation*}
Similarly, 
\[
\Lambda\theta^{\epsilon}_{t} \in L^{4/3}\left(0,T; W^{-2,1}(\mathbb{T})\right).
\]
Moreover, by Sobolev embedding
\[
L^{2}\left(0,T; H^{1}(\mathbb{T})\right)  \subset L^{2}\left(0,T; C^{\alpha}(\mathbb{T})\right), \quad 0<\alpha<1/2,
\]
we have
\begin{equation*} 
 \begin{split}
 \left\|\theta^{\epsilon}(\cdot)-\theta^{\epsilon}(\cdot +y)\right\|_{L^{2}(0,T; L^{2}(\mathbb{T}))} \leq C(\theta_0)|y|^{\alpha }\rightarrow 0 \ \text{as $|y| \rightarrow 0$.}
  \end{split}
\end{equation*}
Similarly, 
\[
\left\|\mathcal{H}\theta^{\epsilon}(\cdot)-\mathcal{H}\theta^{\epsilon}(\cdot +y)\right\|_{L^{2}(0,T; L^{2}(\mathbb{T}))} \rightarrow 0 \ \text{as} \ |y| \rightarrow 0.
\]
Therefore, Lemma \ref{lemma:2.2} with
\[
q_i=p_i=2, \quad g^\epsilon=\theta^\epsilon_x, \ \Lambda\theta^\epsilon, \quad h^\epsilon=\mathcal{H}\theta^\epsilon, \ \theta^\epsilon
\]
allows to pass to the limit in 
\[
\int^{\infty}_{0}\int_{\mathbb{T}} \left[\theta^{\epsilon}_{t} +\left(\mathcal{H}\theta^{\epsilon}\right)\theta^{\epsilon}_{x} +\Lambda^{\gamma}\theta^{\epsilon}+\delta \theta^{\epsilon} \Lambda \theta^{\epsilon} \right]\psi dxdt= \int_{\mathbb{T}}\theta^{\epsilon}_{0}(x)\psi(x,0)dx
\]
to obtain a weak solution in the space stated in (\ref{eq:6.6}).

To show the uniqueness of a solution, let $\theta=\theta_{1}-\theta_{2}$. Then, $\theta$ satisfies
\eqn \label{eq:6.7}
\theta_{t}+\nu\Lambda^{\gamma}\theta=-\left(\mathcal{H}\theta\right) \theta_{1x}- \left(\mathcal{H}\theta_{2}\right) \theta_{x}-\delta \theta \Lambda \theta_{1}-\delta \theta_{2}\Lambda \theta,  \quad \theta(0,x)=0.  
\een
We multiply $\theta$ to (\ref{eq:6.7}) and integrate over $\mathbb{T}$. Then, 
\begin{equation*} 
 \begin{split}
 \frac{1}{2}\frac{d}{dt}\left\|\theta\right\|^{2}_{L^{2}}+\nu\left\|\Lambda^{\frac{\gamma}{2}}\theta \right\|^{2}_{L^{2}}  &= \int_{\mathbb{T}} \left[-\left(\mathcal{H}\theta\right) \theta_{1x}- \left(\mathcal{H}\theta_{2}\right) \theta_{x}-\delta \theta \Lambda \theta_{1}-\delta \theta_{2}\Lambda \theta\right]\theta dx \\
 &=\text{I+II+III+IV}.
 \end{split}
\end{equation*}
We first estimate $\text{II +IV}$:
\begin{equation*} 
 \begin{split}
\text{II +IV} & \leq C\left(\left\|\theta_{2}\right\|_{L^{\infty}} +\left\|\mathcal{H}\theta_{2}\right\|_{L^{\infty}} \right)\|\theta\|_{L^{2}} \left\|\theta_{x}\right\|_{L^{2}} \\
& \leq C(\nu) \left(\left\|\theta_{2}\right\|^{2}_{L^{\infty}} +\left\|\mathcal{H}\theta_{2}\right\|^{2}_{L^{\infty}} \right)\|\theta\|^{2}_{L^{2}} +\frac{\nu}{4} \left\|\theta_{x}\right\|^{2}_{L^{2}}.
\end{split}
\end{equation*}
To estimate $\text{I}$, we do the integration by parts to obtain
\[
\text{I}\leq C\left\|\theta_{1}\right\|_{L^{\infty}} \|\theta\|_{L^{2}} \left\|\theta_{x}\right\|_{L^{2}} \leq \frac{\nu}{8} \left\|\theta_{x}\right\|^{2}_{L^{2}} + C(\nu) \left\|\theta_{1}\right\|^{2}_{L^{\infty}} \|\theta\|^{2}_{L^{2}}.
\]
To estimate $\text{III}$, we use $\Lambda \theta_{1}=(\mathcal{H}\theta_{1})_{x}$ and do the integration by parts to obtain
\[
\text{III}\leq C\left\|\mathcal{H}\theta_{1}\right\|_{L^{\infty}} \|\theta\|_{L^{2}} \left\|\theta_{x}\right\|_{L^{2}} \leq \frac{\nu}{8} \left\|\theta_{x}\right\|^{2}_{L^{2}} + C(\nu) \left\|\mathcal{H}\theta_{1}\right\|^{2}_{L^{\infty}} \|\theta\|^{2}_{L^{2}}.
\]
Since 
\[
\|\theta_{x}\|_{L^{2}}\leq \left\|\Lambda^{\gamma/2}\theta\right\|_{L^{2}}, \quad \left\|\theta_{i}\right\|_{L^{\infty}} + \left\|\mathcal{H}\theta_{i}\right\|_{L^{\infty}}  \leq C\left\|\theta_{i}\right\|_{H^{1}}, \ \text{for $i=1,2$},
\]
we conclude that 
\eqn \label{eq:6.8}
 \frac{d}{dt}\left\|\theta\right\|^{2}_{L^{2}} \leq C(\nu) \left(\|\theta_{1}\|^{2}_{H^{1}} +  \|\theta_{2}\|^{2}_{H^{1}}\right)\|\theta\|^{2}_{L^{2}}
\een
which implies that $\theta=0$ in $L^{2}$.

\section{Proof of Theorem \ref{thm:3.4} and Corollary \ref{cor:3.5}} \label{sec:7}

\subsection*{Proof of Theorem \ref{thm:3.4}}
We consider the equation (\ref{nd dqg}):
\[
\theta_{t} +(1-\delta)u\cdot\nabla\theta + \delta \nabla\cdot(\theta \mathcal{R}\theta)=0.
\]
As the equation (\ref{1d qg}), the entropy \eqref{entropy} satisfies that 
\[
\frac{d}{dt}\mathfrak{E}(\theta)=\int_{\mathbb{T}^n} \theta_{t} \log \theta dx=\int_{\mathbb{T}^n} ((1-\delta)u+\delta\mathcal{R}\theta)\cdot \nabla\theta dx=-\delta\left\|\Lambda^{1/2}\theta\right\|_{L^2}^2.
\]
Therefore, we have $\theta  \in \mathcal{A}$. We now construct a sequence of solutions $\left(\theta^{\epsilon}\right)$ by solving
\[
\theta^\epsilon_{t} +(1-\delta)u^\epsilon\cdot\nabla\theta^\epsilon + \delta \nabla\cdot(\theta^\epsilon \mathcal{R}\theta^\epsilon)=\epsilon\Delta \theta^\epsilon,\quad \theta^\epsilon_{0}=\rho_\epsilon \ast \theta_0.
\]
Then, $(\theta^{\epsilon})$ satisfies that
\eqn \label{eq:7.1}
\frac{d}{dt}\mathfrak{E}(\theta^{\epsilon}) +\left\|\Lambda^{1/2}\theta^{\epsilon}\right\|^{2}_{L^{2}} + 4\epsilon\int_{\mathbb{T}^{n}}\left|\nabla \left(\sqrt{\theta^{\epsilon}}\right)\right|^{2}dx= 0.
\een
Integrating (\ref{eq:7.1}) in time,
\eqn\label{eq:7.2}
\mathfrak{E}(\theta^{\epsilon}(t))  +\int^{t}_{0}\left\|\Lambda^{1/2}\theta^{\epsilon}(s)\right\|^{2}_{L^{2}}ds  + 4\epsilon \int^{t}_{0}\int_{\mathbb{T}^{n}}\left|\nabla \left(\sqrt{\theta^{\epsilon}}(s)\right)\right|^{2}dxds \leq \mathfrak{E}(\theta^{\epsilon}_{0}).
\een
Since $x\log x -x+1 \leq x^{s+1}+1$ for $x\ge 0$, we can bound the last term in (\ref{eq:7.2}) by
\[
\mathfrak{E}(\theta^{\epsilon}_{0}) \leq  (2\pi)^{n}+ \left\| \theta^{\epsilon}_{0}\right\|^{s+1}_{L^{s+1}}. 
\]
Therefore, the sequence $\left(\theta^{\epsilon}\right)$ is uniformly bounded in $\mathcal{A}_T$. Using this bound, we first treat the two dimensional case. From Sobolev embedding, we have uniform bounds 
\[
\mathcal{R}\theta^\epsilon, \ u^\epsilon, \ \theta^\epsilon\in L^2\left(0,T; L^4(\mathbb{T}^2)\right).
\] 
Moreover, by interpolating $L^\infty\left(0,T; L^1(\mathbb{T}^2)\right)$ and $L^2\left(0,T; L^4(\mathbb{T}^2)\right)$, we  have uniform bounds 
\[
\mathcal{R}\theta^\epsilon, \ u^\epsilon, \ \theta^\epsilon\in L^6\left(0,T; L^{4/3}(\mathbb{T}^2)\right).
\] 
This implies that 
\[
u^{\epsilon}\theta^{\epsilon}, \ \theta^{\epsilon}\mathcal{R}\theta^{\epsilon} \in L^{3/2}\left(0,T; L^1(\mathbb{T}^2)\right)
\]
and thus
\[
\theta^\epsilon_t\in L^{3/2} \left(0,T; H^{-2}(\mathbb{T}^2)\right).
\]
Using Lemma \ref{lemma:2.1} with 
\[
X_0=L^2 \left(0,T; H^{1/2}(\mathbb{T}^2)\right),\quad X=L^2\left(0,T; L^2(\mathbb{T}^2)\right),\quad X_1=L^{3/2}\left(0,T; H^{-2}(\mathbb{T}^2)\right)
\]
we can pass to the limit in the weak formulation to obtain a weak solution in $\mathcal{A}_{T}$. Similarly,  in three dimensions, we have
\[
u^{\epsilon}\theta^{\epsilon}, \ \theta^{\epsilon}\mathcal{R}\theta^{\epsilon} \in L^{4/3}\left(0,T; L^1(\mathbb{T}^3)\right), \quad \theta^\epsilon_t\in L^{4/3} \left(0,T; H^{-2}(\mathbb{T}^3)\right).
\]
Using Lemma \ref{lemma:2.1} with
\[
X_0=L^2\left(0,T; H^{1/2}(\mathbb{T}^3)\right),\quad X=L^2\left(0,T; L^2(\mathbb{T}^3)\right),\quad X_1=L^{4/3}2\left(0,T; H^{-2}(\mathbb{T}^3)\right),
\]
we complete the proof of Theorem \ref{thm:3.4}.

\subsection*{Proof of Corollary \ref{cor:3.5}}
We notice that the hypothesis $m(0)=0$ together with  $\theta^\epsilon\in  L^{2} \left((0,T); H^{1/2}(\mathbb{T}^{n})\right)$ implies $u^\epsilon\in L^{2} \left((0,T); H^{1/2+s}(\mathbb{T}^{n})\right)$ for some $s>0$. This is enough to follow the argument in the proof of Theorem \ref{thm:3.4} to complete the proof of Corollary \ref{cor:3.5}.

\section{Proof of Theorem \ref{thm:3.6}}\label{sec:8}

\noindent
We finally consider the equation 
\[
\theta_{t}+u\cdot \nabla \theta +\nu\Lambda^{\gamma} \theta=0
\]
with the following entropy
\[
\mathfrak{E}(\theta)=\int_{\mathbb{T}^{n}}(\theta+1)\log(\theta+1)dx.
\]
Since $\theta\ge 0$, we have $\mathfrak{E}(\theta)\ge 0$. Let's start with the \emph{a priori} estimates. The direct computation yields that 
\eqn \label{eq:8.1}
\frac{d}{dt}\mathfrak{E}(\theta)=-\nu \int_{\mathbb{T}^{n}} \left(\Lambda^{\gamma}\theta\right) \log(\theta+1)dx.
\een
To estimate the right hand side of (\ref{eq:8.1}) we symmetrize the integral using the representation of $\Lambda^{\gamma}$ (\ref{lambda gamma}):
\begin{equation*} 
 \begin{split}
 \frac{d}{dt}\mathfrak{E}(\theta)& + \frac{c_{\gamma,n}\nu}{2}\sum_{k\in \mathbb{Z}^{n}} \int_{\mathbb{T}^{n}}\int_{\mathbb{T}^{n}}\frac{\theta(x)-\theta(y)}{|x-y+2k\pi|^{n+\gamma}}\log \left[\frac{\theta(x)+1}{\theta(y)+1}\right]dydx =0. 
  \end{split}
\end{equation*}
Since  $(X-Y)(\log X-\log Y) \ge C(\log X-\log Y)^{2}$ for $X\ge 1$ and $Y\ge 1$, we obtain
\[
\frac{d}{dt}\mathfrak{E}(\theta)+ C(\nu, \gamma, n)\left\| \Lambda^{\gamma/2}\log (\theta+1)\right\|^{2}_{L^{2}}\leq 0. 
\]
To obtain the diffusion, we compute
\begin{equation*} 
 \begin{split}
 &c(\gamma,n,\nu)\sum_{k\in \mathbb{Z}^{n}} \int_{\mathbb{T}^{n}}\int_{\mathbb{T}^{n}}\frac{(\theta(x)-\theta(y))^2}{|x-y+2k\pi|^{n+\gamma}}\frac{(\log(1+\theta(x))-\log(1+\theta(y)))^2}{(\theta(x)-\theta(y))^2}dydx\\
&\geq \frac{c(\gamma,n,\nu) }{1+\left\| \theta_{0}\right\|^{2}_{L^{\infty}}}\|\Lambda^{\gamma/2}\theta\|_{L^2}^2. 
  \end{split}
\end{equation*}
This implies that 
\[
\mathfrak{E}(\theta(t)) + \frac{c(\gamma,n,\nu) }{1+\left\| \theta_{0}\right\|^{2}_{L^{\infty}}}\|\Lambda^{\gamma/2}\theta\|_{L^2}^2\leq \mathfrak{E}(\theta_{0}) \leq (2\pi)^{n}+ \left\| \theta_{0}\right\|^{s+1}_{L^{s+1}}
\]
for all $t>0$ and thus we conclude $\theta \in \mathcal{B}_{T}$ from (\ref{eq:8.1}). We now construct a sequence of solutions $\left(\theta^{\epsilon}\right)$ by solving
\[
\theta^\epsilon_{t} +u^\epsilon\cdot\nabla\theta^\epsilon + \nu \Lambda^{\gamma}\theta^{\epsilon}=\epsilon\Delta \theta^\epsilon,\quad \theta^\epsilon_{0}=\rho_\epsilon \ast \theta_0.
\]

For such a solution, the following inequality holds:
\[
\|\theta^\epsilon(t)\|_{L^1\left(\mathbb{T}^{n}\right)}=\|\theta^\epsilon_0\|_{L^1\left(\mathbb{T}^{n}\right)}\leq \|\theta_0\|_{L^1\left(\mathbb{T}^{n}\right)}.
\]
Since $\theta^\epsilon(t)$ is smooth (in space and time), the function $\|\theta^\epsilon(t)\|_{L^\infty\left(\mathbb{T}^{n}\right)}=\theta^\epsilon(x_t)$ is Lipschitz. Using Rademacher Theorem, we conclude that $\theta(x_t)$ is differentiable almost everywhere. For each $t$, let $x_t$ be the point of maximum. Then, 
\[
\frac{d}{dt}\|\theta^\epsilon(t)\|_{L^\infty\left(\mathbb{T}^{n}\right)}=\pat \theta^\epsilon(x_t).
\]
We now estimate nonlocal terms. We take  a positive number $0<r<\pi$ and define 
\[
\mathcal{U}_1=\left\{\eta\in[-r,r]^n: \theta^\epsilon(x_t)-\theta^\epsilon(x_t-\eta)>\theta^\epsilon(x_t)/2 \right\},
\]
and $\mathcal{U}_2=[-r,r]^n\setminus \mathcal{U}_1$. Then, we have 
\[
\|\theta_0\|_{L^1\left(\mathbb{T}^{n}\right)}\geq\|\theta^\epsilon_0\|_{L^1\left(\mathbb{T}^{n}\right)}=\int_{\mathbb{T}^{n}} \theta^\epsilon(x_t-\eta)d\eta\geq\int_{\mathcal{U}_2} \theta^\epsilon(x_t-\eta)d\eta\geq \frac{\theta^\epsilon(x_t)}{2} \left|\mathcal{U}_2\right|
\]
or equivalently, 
\[
(2r)^n-\frac{2\|\theta_0\|_{L^1\left(\mathbb{T}^{n}\right)}}{\theta^\epsilon(x_t)}\leq (2r)^n- \left|\mathcal{U}_2\right|= \left|\mathcal{U}_1\right|.
\]
This implies that 
\begin{equation*} 
 \begin{split}
 & \frac{1}{c_{\gamma,n}}\Lambda^\gamma \theta^\epsilon(x_t)\geq \text{P.V.}\int_{\mathbb{T}^{n}}\frac{\theta^\epsilon(x_t)-\theta^\epsilon(x_t-y)}{|y|^{n+\gamma}}dy \geq\int_{\mathcal{U}_1}\frac{\theta^\epsilon(x_t)-\theta^\epsilon(x_t-y)}{|y|^{n+\gamma}}dy\\
 & \geq\frac{\theta^\epsilon(x_t)}{2r^{n+\gamma}} \left|\mathcal{U}_1\right| \geq \frac{\theta^\epsilon(x_t)}{r^{n+\gamma}}\left(2^{n-1}r^n-\frac{\|\theta_0\|_{L^1\left(\mathbb{T}^{n}\right)}}{\theta^\epsilon(x_t)}\right).
  \end{split}
\end{equation*}
We now choose $r$ as follows: 
\[
r=\left(\frac{1}{2^{n-2}}\frac{\|\theta_0\|_{L^1\left(\mathbb{T}^{n}\right)}}{\theta^\epsilon(x_t)}\right)^{1/n}.
\]
We assume $r\leq\pi$ for the moment. In this case, we obtain that
\[
\frac{1}{c_{\gamma,n}}\Lambda^{\gamma} \theta^\epsilon(x_t)\geq \frac{\|\theta_0\|_{L^1\left(\mathbb{T}^{n}\right)}}{\left(\frac{1}{2^{n-2}}\frac{\|\theta_0\|_{L^1\left(\mathbb{T}^{n}\right)}}{\theta^\epsilon(x_t)}\right)^{1+\gamma/n}}.
\]
This bound implies that
\[
\frac{d}{dt}\|\theta^\epsilon(t)\|_{L^\infty\left(\mathbb{T}^{n}\right)}\leq -\nu c_{\gamma,n}\Lambda^\gamma \theta^\epsilon(x_t) +\epsilon\Delta\theta^\epsilon  \leq -\nu c_{\gamma,n}\frac{\|\theta^\epsilon(t)\|_{L^\infty\left(\mathbb{T}^{n}\right)}^{1+\gamma/n}}{\|\theta_0\|^{\gamma/n}_{L^1\left(\mathbb{T}^{n}\right)}},
\]
or equivalently,
\[
\|\theta^\epsilon(t)\|_{L^\infty(T^n)}\leq \left(\frac{\nu c_{\gamma,n}}{\|\theta_0\|^{\gamma/n}_{L^1\left(\mathbb{T}^{n}\right)}}\gamma t+\frac{1}{\|\theta^\epsilon_0\|_{L^\infty}^\gamma}\right)^{-1/\gamma}\\
\leq \frac{\|\theta_0\|^{1/n}_{L^1\left(\mathbb{T}^{n}\right)}}{\left(C(\nu, \gamma,n) \right)^{1/\gamma}}t^{-1/\gamma}.
\]
If $r$ is bigger than $\pi$, 
\[
\frac{1}{2^{n-2}\pi^n}\|\theta_0\|_{L^1\left(\mathbb{T}^{n}\right)}>\|\theta^\epsilon (t)\|_{L^\infty\left(\mathbb{T}^{n}\right)}.
\]
As a consequence, for $t>0$, we have
\[
\|\theta^\epsilon(t)\|_{L^\infty\left(\mathbb{T}^{n}\right)}\leq \max\left\{\frac{\|\theta_0\|^{1/n}_{L^1\left(\mathbb{T}^{n}\right)}}{\left(C(\nu, \gamma,n) \right)^{1/\gamma}}t^{-1/\gamma}, \ \frac{1}{2^{n-2}\pi^n}\|\theta_0\|_{L^1\left(\mathbb{T}^{n}\right)}\right\}.
\]
In particular, for $t\ge \tau$, 
\eqn \label{eq:9.1}
\|\theta^\epsilon(t)\|_{L^\infty\left(\mathbb{T}^{n}\right)}\leq \max\left\{\frac{\|\theta_0\|^{1/n}_{L^1\left(\mathbb{T}^{n}\right)}}{\left(\nu C_{\gamma,n}\gamma \right)^{1/\gamma}}\tau^{-1/\gamma},\frac{1}{2^{n-2}\pi^n}\|\theta_0\|_{L^1\left(\mathbb{T}^{n}\right)}\right\}=:C(\tau, \nu, \gamma,n,\gamma).
\een

Then, $(\theta^{\epsilon})$ is uniformly bounded in $\mathcal{B}_{T}$. Following the proof of Theorem \ref{thm:3.1}, we obtain a weak solution $\theta$ in $\mathcal{B}_{T}$. Since 
\begin{equation*} 
 \begin{split}
 \theta^{\epsilon}_{t} =-u^\epsilon\cdot\nabla\theta^\epsilon - \nu \Lambda^{\gamma}\theta^{\epsilon}+ \epsilon\Delta \theta^\epsilon \in L^{2}\left(0,T; H^{-2}(\mathbb{T}^{n})\right)
  \end{split}
\end{equation*}
uniformly in $\epsilon>0$, we have
\[
\theta^{\epsilon} \in C\left(0,T; H^{-2}(\mathbb{T})\right).
\]
Therefore, we recover $\theta_{0}$ in $H^{-2}(\mathbb{T}^{n})$.

\section{Appendix}

\subsection{Lyapunov functions}
For the equation (\ref{1d qg}), we have two additional Lyapunov functions defined in terms of 
\eqn
-\Lambda^{-1}\theta =\frac{1}{\pi}\int_{\mathbb{T}} \log \left|\sin\left(\frac{x-y}{2}\right)\right|\theta(y)dy.
\een
When $\theta_{0}\ge 0$, by the minimum principle, $-\Lambda^{-1}\theta \leq 0$. The first Lyapunov function is 
\eqn
\mathcal{L}_{1}(\theta)=\int_{\mathbb{T}} \left[\theta \left(\log \left|\Lambda^{-1}\theta\right|\right)+M\right]dx, \quad M=\|\theta_{0}\|_{L^{\infty}} \left\|\log \left|\Lambda^{-1}\theta_{0}\right|\right\|_{L^{\infty}},
\een
where we add $M$ to make $\mathcal{L}_{1}(\theta)$ to be non-negative. We show that $\mathcal{L}_{1}(\theta)$ is decreasing in time using the fact $-\Lambda^{-1}f_{x}=\mathcal{H}f$: 
\begin{equation*} 
 \begin{split}
 & \frac{d}{dt}\mathcal{L}_{1}(\theta)= \int_{\mathbb{T}} \left[\theta_{t} \log \left|\Lambda^{-1}\theta\right|\right]dx + \int_{\mathbb{T}} \left[\theta \left(\log \left|\Lambda^{-1}\theta\right|\right)_{t}\right]dx\\
 & =\int_{\mathbb{T}} \left[\left(\theta \mathcal{H}\theta\right) \left(\log \left|\Lambda^{-1}\theta\right|\right)_{x}\right]dx + \int_{\mathbb{T}} \left[\frac{\theta\Lambda^{-1}\theta_{t}}{\Lambda^{-1}\theta}\right]dx =\int_{\mathbb{T}} \left[\frac{-\theta\left(\mathcal{H}\theta\right)^{2}}{\Lambda^{-1}\theta}+\frac{\theta \mathcal{H} \left(\theta\mathcal{H}\theta\right)}{\Lambda^{-1}\theta}\right]dx\\
 & =-\frac{1}{2}\int_{\mathbb{T}} \left[\frac{\theta}{\Lambda^{-1}\theta} \left((\mathcal{H}\theta)^{2}+(\theta)^{2}+\langle \theta \rangle^{2}\right)\right]dx\leq 0.
  \end{split}
\end{equation*}
The second Lyapunov function is 
\eqn
\mathcal{L}_{2}(\theta)=\int_{\mathbb{T}} \left[\theta e^{\Lambda^{-1}\theta}\right]dx.
\een
We show that $\mathcal{L}_{2}(\theta)$ is exponentially decreasing in time:
\begin{equation*} 
 \begin{split}
 & \frac{d}{dt}\mathcal{L}_{2}(\theta)= \int_{\mathbb{T}} \left[\theta_{t} e^{\Lambda^{-1}\theta}\right]dx + \int_{\mathbb{T}} \left[\theta \left(e^{\Lambda^{-1}\theta}\right)_{t}\right]dx\\
 & =\int_{\mathbb{T}} \left[\left(\theta \mathcal{H}\theta\right) \left( e^{\Lambda^{-1}\theta}\right)_{x}\right]dx + \int_{\mathbb{T}} \left[\theta e^{\Lambda^{-1}\theta} \mathcal{H} \left(\theta \mathcal{H}\theta\right)\right]dx\\
 & =-\int_{\mathbb{T}} \left[\theta\left(\mathcal{H}\theta\right)^{2} e^{\Lambda^{-1}\theta}\right]dx + \frac{1}{2}\int_{\mathbb{T}} \left[\theta e^{\Lambda^{-1}\theta}\left((\mathcal{H}\theta)^{2}-\theta^{2}-\langle \theta \rangle^{2}\right)\right]dx\\
 & =-\frac{1}{2}\int_{\mathbb{T}} \left[\theta e^{\Lambda^{-1}\theta} \left((\mathcal{H}\theta)^{2}+\theta^{2}\right)\right]dx -\frac{\langle \theta_{0} \rangle^{2}}{2} \mathcal{L},
  \end{split}
\end{equation*}
where we use $\langle \theta \rangle=\langle \theta_{0} \rangle$. Therefore, 
\eqn
\mathcal{L}_{2}(\theta)\leq e^{-\frac{\langle \theta_{0} \rangle^{2}}{2}t} \int_{\mathbb{T}} \left[\theta_{0} e^{\Lambda^{-1}\theta_{0}}\right]dx.
\een 
We note that the bound of $\mathcal{L}_{2}$ implies the bound of $\theta$ in $\dot{H}^{-1/2}(\mathbb{T})$: 
\[
\|\theta(t)\|_{\dot{H}^{-1/2}(\mathbb{T})}\leq \int_{\mathbb{T}} \theta \left(1+\Lambda^{-1}\theta\right)dx\leq  \mathcal{L}_{2}(\theta)<C(\theta_{0}).
\]

\section*{Acknowledgments}
H.B. and RGB acknowledge Prof. H. Dong for his valuable comments that highly improve the manuscript. The authors also gratefully acknowledge the support by the Department of Mathematics at UC Davis where this research was performed. The second author is partially supported by the grant MTM2011-26696 from the former Ministerio de Ciencia e Innovaci\'on (MICINN, Spain).

\bibliographystyle{abbrv}

\end{document}